\documentclass[11pt,reqno,a4paper]{amsart}
\usepackage{euscript,verbatim}

\newtheorem{theorem}{Theorem}[section]
\newtheorem{proposition}[theorem]{Proposition}

\newtheorem{corollary}[theorem]{Corollary}

\newtheorem{lemma}[theorem]{Lemma}

\newcommand{\beha}{\begin{enumerate}}
\newcommand{\behe}{\end{enumerate}}
\renewcommand{\epsilon}{\varepsilon}

\def\htop{h_{\text{top}}}
\def\inte{{\rm int}\ }
\def\vol{{\rm vol} }
\newcommand{\bR}{{\mathbb R}}
\newcommand{\bC}{{\mathbb C}}
\newcommand{\bP}{{\mathbb P}}

\newcommand{\bN}{{\mathbb N}}

\newcommand{\cM}{\EuScript{M}}

\newcommand{\RE}[1]{(\ref{#1})}

\DeclareMathSymbol{\varnothing}{\mathord}{AMSb}{"3F}
\renewcommand{\emptyset}{\varnothing}

\begin{document}

\title[Polynomial Automorphisms of~$\bC^n$]{Filtrations, Hyperbolicity and Dimension for polynomial Automorphisms of~$\bC^n$}

\begin{abstract}
In this paper we study the dynamics of  regular polynomial
automorphisms of $\bC^n$. These maps provide a natural
generalization of complex H\'enon maps in $\bC^2$ to higher
dimensions. For a given  regular polynomial automorphism $f$
we construct a filtration in $\bC^n$ which has particular
escape properties for the orbits of $f$. In the case when $f$ is
hyperbolic we obtain a complete description of its orbits.
In the second part of the paper we study
the Hausdorff and box dimension of the Julia sets of $f$. We show
that the Julia set $J$  has positive box dimension, and (provided
$f$ is not volume preserving) that the filled-in Julia set $K$ has box
dimension strictly less than $2n$. Moreover, if $f$ is
 hyperbolic, then the Hausdorff dimension of the forward/backward
Julia set $J^\pm$  is strictly less than $2n$.
\end{abstract}

\date{\today}

\thanks{C.W. was partially supported by the Center for Mathematical
Analysis, Geometry, and Dynamical Systems, through FCT's Funding
Program.}

\author{Rasul Shafikov} \address{Department of Mathematics, SUNY,
Stony Brook, NY 11794, USA}
\email{shafikov@math.sunysb.edu}

\author{Christian Wolf}\address{Departamento de Matem\'atica, Instituto Superior
T\'ecnico, 1049-001 Lisboa, Portugal}
\email{cwolf@math.ist.utl.pt}

%\subjclass{Primary: 11K55, 37C45.}
\maketitle

%-------------------------------------
\section{Introduction}

\bigskip

Let $f$ be a polynomial automorphism of $\bC^n$. We denote by
$\hat f$ the  natural extension of $f$ to a meromorphic map in
$\bP^n$. Let $I^+$ denote the indeterminacy  set of $\hat f$.
Analogously we denote by $I^-$ the indeterminacy set of
$\widehat{f^{-1}}$. We say that  $f$  is regular if $f$ has
degree greater than one and $I^+\cap I^-=\emptyset$.

In the case  $n=2$  the class of regular automorphisms consists of
polynomial automorphisms with nontrivial dynamics, i.e.,
finite compositions of generalized H\'enon maps (see for instance
\cite{BS},\cite{FM},\cite{FS}). In fact, regular polynomial
automorphisms can be considered as a natural generalization of
complex H\'enon maps to higher dimensions. Higher dimensional
regular maps are for instance the so-called shift-like
automorphisms studied by Bedford and Pambuccian \cite{BP}. For
further examples we refer to Section 2. We point out that unlike
the two-dimensional case, for $n>2$ there exist polynomial
automorphisms with nontrivial dynamics which are not regular (see
for instance \cite{CF}).

The notion of regular polynomial automorphisms was introduced by
Sibony \cite{Si}, who comprehensively studied  these maps, in
particular by using 
methods from pluripotential theory.

In this paper we study the dynamics of regular polynomial
automorphisms from a different point of view: We introduce the
notion of hyperbolicity for a regular polynomial automorphism $f$
and study its dynamics. In particular, we classify the orbits of
$f$ analogously to the case of complex H\'enon maps in \cite{BS}. 
Finally, we study the Hausdorff and box dimension of
the Julia sets of $f$. We derive estimates for these dimensions in
the hyperbolic as well as in the non-hyperbolic case.

We will now describe our results in more detail:

Let $f$ be a regular polynomial automorphism of $\bC^n$. We define
$K^{\pm}= \{ p\in \bC^n: \{f^{\pm k} (p): k\in \bN\} {\rm\ is \ bounded}\}$,
and the filled-in Julia set  by $K=K^+\cap K^-$. Furthermore, we
define the sets $J^\pm=\partial K^\pm$ and $J= J^+\cap J^-$. The
set $J^\pm$ is called the forward/backward Julia set and $J$ is
the Julia set of $f$ (see Section 2 for details).

We construct a filtration of $\bC^n$ which has particular escape
properties for the orbits of $f$ (see Proposition~\ref{l1.5}). For
$n=2$ the existence of  a filtration was already shown  by Bedford
and Smillie \cite{BS}. In this case it proved to be a useful tool
for analyzing the dynamics of $f$.

We apply our filtration to study hyperbolic maps. We say that $f$
is hyperbolic if its Julia set $J$ is a hyperbolic set and the
periodic points are dense in $J$. It turns out  that hyperbolicity
implies that $f$ is Axiom~A (see Corollary~\ref{coraxioma}). The
latter is the classical notion for hyperbolic diffeomorphisms. We
obtain a complete description for the possibilities of the orbits
in the case of a hyperbolic map $f$. The following result is a
consequence of Proposition \ref{l:hyper} and Theorems \ref{k+} and
\ref{thghj}.
\begin{theorem}
Let $f$ be a hyperbolic regular polynomial automorphism of $\bC^n$
and let $p$ be a point in $\bC^n$. Then one of the following
exclusive properties holds.
\begin{enumerate}\label{thmainhyp}
\item[(i)] There exists $q\in J$ such that $|f^k(p)-f^k(q)|\to 0$
as $k\to \infty$;
\item[(ii)]
 There exists an attracting periodic point $\alpha$ of $f$ such that $|f^k(p)-f^k(\alpha)|\to 0$
 as
$k\to \infty$;
\item[(iii)]
$\{f^k(p): k\in\bN\}$ converges to $\infty$ as $k\to\infty$.
\end{enumerate}
\end{theorem}

The inverse of $f$ is also a regular polynomial automorphism;
therefore, Theorem \ref{thmainhyp} also holds for $f^{-1}$.
However, since $f$ has constant jacobian, attracting periodic
points can exist only either for $f$ or for $f^{-1}$. Theorem
\ref{thmainhyp} implies that in order to understand the
``complicated" dynamics  of a hyperbolic map, it is sufficient to
understand the dynamics on its Julia set.

In the second part of this paper we derive estimates for the Hausdorff
and  box dimension of the Julia sets. Let $d$ be the degree of
$f$; then we denote by $l-1$ the dimension of $I^-$ (as an
algebraic variety). We show that $J$ carries the full entropy of
$f$, that is, $\htop(f|J)=l\log d$ (see Theorem \ref{prohtop}),
where $\htop$ denotes the topological entropy. As a consequence,
the upper box dimension of the Julia set is strictly positive (see
Corollary \ref{corkk3}). On the other hand, if $f$ is
not volume preserving, then the upper box dimension of $K$ is
strictly smaller than $2n$ (see Corollary \ref{corkk4}).

It is a widely studied problem in one-dimensional complex dynamics
to determine whether the Hausdorff dimension of the Julia set of a
rational map is strictly less than two (see for instance \cite{U}
for an overview). We solve the analogous problem in the case of
hyperbolic regular polynomial automorphisms of $\bC^n$, $n\ge 2$; 
namely we show that the Hausdorff dimension of $J^\pm$ is strictly 
less than $2n$. More precisely, we derive an upper bound for the 
Hausdorff dimension of $J^\pm$ which is given in terms of topological
pressure (see Theorem \ref{jpmk4}). This upper bound is strictly
smaller than $2n$. Our Theorem improves a result of Bowen
\cite{Bo}, which says that $J^\pm$ has zero Lebesgue measure. It
should be noted that for $n=2$ it is possible to construct
hyperbolic maps whose forward/backward Julia sets have Hausdorff
dimension arbitrarily close to $4=2n$ (see \cite{W3}).

This paper
is organized as follows. In Section 2 we present basic facts
about regular polynomial automorphisms of $\bC^n$. In Section 3 we
construct a filtration with particular escape properties for
the orbits. Section 4 is devoted to the analysis of hyperbolic
maps. Finally, we study in Section 5 the Hausdorff and the box
dimension of the Julia sets.

For $n=2$ it is shown in \cite{BS} that hyperbolicity of $J$
already implies that $f$ is a hyperbolic map, i.e., $J$ being
hyperbolic implies the density of the periodic points in $J$. In
particular,  this provides a weaker definition of a hyperbolic map
in the case $n=2$. It would be interesting to know
whether the analogous result
holds for $n>2$, see also the remark at the end of Section 4.

%% ----------------------- section ---------------------------

\section{Regular polynomial automorphisms}

In this section we give an  introduction to the dynamics of regular
polynomial automorphisms of $\bC^n$. This class of maps was
studied in detail by Sibony in \cite{Si}. We refer to this article
for the proofs of the results presented below.

Let $f$ be a polynomial automorphism of $\bC^n, n\geq 2$. Then $f$
admits an extension to a meromorphic map $\hat f: \bP^n\to\bP^n$.
Let $\pi: \bC^{n+1}\to\bP^n$ be the canonical projection and let
$F$ be the homogeneous polynomial map in $\bC^{n+1}$ corresponding
to $\hat f$, i.e., $\hat f=\pi\circ F\circ \pi^{-1}$. Then
$I^+=\pi\circ F^{-1}(0)$ is the indeterminacy set of $\hat f$.
Analogously, $I^-$ denotes the indeterminacy set of $\widehat
{f^{-1}}$. The sets $I^+$ and $I^-$ are algebraic varieties in
$\bP^n$ of codimension at least two which are contained in the
hypersurface at infinity, which is denoted by $H_0$.

We write $f=(f_1,\dots,f_n)$. Let  $\deg f_i$ denote the polynomial
degree of $f_i$. Then $d=\deg f=\max\{\deg f_1,\dots,\deg f_n\}$ is
the degree of $f$. A polynomial automorphism $f$ of $\bC^n$ is called {\it regular} if $d>1$ and $I^+\cap I^- =\varnothing$.

Throughout this paper $f$ will be a regular polynomial
automorphism of $\bC^n$. Note that the complex Jacobian  $\det\
Df$ is constant in $\bC^n$. Therefore, we can restrict our
considerations to the volume decreasing case $(|\det\ Df|<1)$ and
to the volume preserving case $(|\det\ Df|=1)$, as otherwise we
simply consider $f^{-1}$.

Examples of regular polynomial automorphisms are the well known
family of generalized H\'enon maps in $\bC^2$ and shift-like
mappings in $\bC^n$ (see \cite{BP}). Forn\ae ss and Wu \cite{FW}
gave a classification of polynomial automorphisms of $\bC^3$ of
degree at most two. According to this classification automorphisms
with non-trivial dynamics fall into five classes. One can easily
verify that  two of these classes, i.e., the maps of the form
\begin{eqnarray*}
h_1 = (P(x,y)+az, Q(y)+x, y),\\ h_2 = (P(x,y)+az, Q(x)+by,x),
\end{eqnarray*}
are regular, provided that $ab\not= 0$, the degree of $P$ is two in 
each variable, and the degree of $Q$ is also two.

For a regular polynomial automorphism $f$ of $\bC^n$, there exists
an integer $l>0$ such that  $(\deg f)^l = (\deg f^{-1})^{n-l}$.
For example, if $h_1$ and $h_2$ are regular, then $\deg
h_1^{-1}=\deg h_2^{-1} =4$. Moreover, it follows that for every
regular polynomial automorphism we have $\dim I^-=l-1$ and $\dim
I^+=n-l-1$.

Let $K^\pm,K, J^\pm$ and $J$ be defined as in the introduction.
Then all of these  are $f$-invariant sets, and
\begin{equation}\label{eqwerwowas}
\overline{K^{\pm}} = K^{\pm} \cup I^{\pm},
\end{equation}
where $\overline{K^{\pm}}$ denotes the closure in $\bP^n$.
While $K^\pm$ and $J^\pm$ are closed and unbounded in $\bC^n$,
the sets $K$ and $J$ are compact. If $n=2$, then the set $J^\pm$ is
connected (see \cite{BS2}).

Furthermore, the distance between $f^k(p)$ and $I^-$ tends
uniformly to zero on compact subsets of $\bC^n\setminus K^+$.  On
the other hand, the distance between $f^{k}(p)$ and $K$ tends
uniformly to zero on compact subsets of $K^+$, and the family
$\{f^k: k\in\bN\}$ is equicontinuous in the interior of $K^+$,
${\rm int\ } K^+$. The analogous properties hold for $f^{-1}$ and
$K^-$.

%%------------------------ section ---------------------------

\section{Filtration Properties}

In this  section we construct for a regular polynomial
automorphism $f$ a filtration in $\bC^n$ that exhibits particular
escape properties for the orbits of $f$. Our approach is motivated
by the work of Bedford and Smillie \cite[\S 2]{BS} in the case of
generalized H\'enon maps.

\begin{proposition}[\bf Filtration]\label{l1.5}
Let $f$ be a regular polynomial automorphism of $\bC^n$. Then
there exist a compact set $V\subset \bC^n$ with $K\subset\inte V$
and sets $V^+, V^-\subset\bC^n$ with $\bC^n=V\cup V^-\cup V^+$,
such that the following hold.
\begin{enumerate}
\item[(i)] $f(V^-)\subset V^-$;
\item[(ii)] $f(V^-\cup V)\subset V^- \cup V$;
\item[(iii)] $f^{-1}(V^+)\subset V^+$;
\item[(iv)] $f^{-1}(V^+\cup V)\subset V^+ \cup V$.
\end{enumerate}
\end{proposition}

\begin{proof}
Let $V$ be a closed polydisk of sufficiently large radius such
that $K\subset \inte V$. Let $\hat V^+$ and $\hat V^-$ be open
sets in $\bP^n$ for which the following properties hold:
\begin{equation}\label{vpm}
\begin{array}{l}
(a)\  I^+\subset \hat V^+, I^-\subset \hat V^-,\\ (b)\
\widehat{f^{-1}}(\hat V^+)\subset \hat V^+, \ \hat f(\hat
V^-)\subset\hat V^-,\\ (c)\  \hat V^{\pm}\cap  V=\varnothing, \ \
\hat V^{\pm} \cap K^{\mp} =\varnothing.
\end{array}
\end{equation}
We note that property (b) above can be satisfied, since $I^+$ is
an attracting set for $\widehat{f^{- 1}}$, and $I^-$ is an
attracting set for $\hat f$ (see \cite[Prop. 2.2.6]{Si}).
Moreover, it is a consequence of \eqref{eqwerwowas} that we can
achieve the identities on the right of (c).

Let $V_0^{\pm}=\hat V^{\pm}\cap \bC^n$. To construct $V^+$ and
$V^-$ we define
\begin{equation}
V_k^+ = f^{k}(V_{0}^+), \ V_k^- = f^{-k}(V_{0}^-),
\end{equation}
where $k\ge 1$, and at each step, if $V_k^+\cap
V_k^-\not=\emptyset$, we replace the set $V_0^-$ with
$f^{k}(V_k^-\setminus V_k^+)$. To emphasize that $V^{-}_0$ depends
on the number of times this process was applied, we use the
notation $V^{-}_0(k)$, where $k$ is the number of iterations.

We also note that the shrinking of $V^{-}_0(k)$ does not
change the union of $V_k^+$ and $V_k^-$, and  that $V_k^+$ and
$V_k^-$, are disjoint for every $k>0$.
To prove the proposition we will show that there exists $N>0$
such that $\bC^n \setminus V \subset V^{+}_N \cup V^{-}_N $.
We will do this in two steps.

\medskip
\noindent{\it Claim 1. There exists $N>0$ such that for every
$p\in \bC^n \setminus V$ and every $k\ge N$ either $f^k(p)\in \hat
V^-\cap \bC^n$, or $f^{-k}(p)\in \hat V^+\cap \bC^n$. }
\medskip

To prove the claim  we consider a neighborhood $U$ in $\bP^n$ of
the compact set $H_0\setminus(\hat V^-\cup \hat V^+)$.  We recall
that $H_0\subset \bP^n$ denotes the hypersurface at infinity.
Without loss of generality we may assume that $U\cap
I^{\pm}=\varnothing$, and $U\cap V=\varnothing$. Furthermore,
since $\overline{K^+}\cap H_0 = I^+$, we may choose $U$ such that
$\overline{U\setminus \hat V^+}\cap K^+=\varnothing$. From $\hat f
(H_0\setminus I^+)=I^-$ and the uniform convergence of $f^k(p)$ to
$I^-$ on compact subsets of $\bC^n\setminus K^+$ we conclude that
there exists $N_1>0$ such that $\widehat{f^k} (U)\subset \hat V^-$
for $k\ge N_1$. In view of  property (b) in (\ref{vpm}) it follows
that Claim 1 holds for any point $p\in (U\cup \hat V^+ \cup \hat
V^-)\cap \bC^n$ and all $N \ge N_1$. We note, that $H_0\subset U\cup \hat V^+
\cup \hat V^-$.

We consider now the set $D = \bC^n \setminus (U\cup \hat V^+ \cup
\hat V^- \cup V)$. Since $K\subset \inte V$, $K^\pm$ is closed and
in view of \eqref{eqwerwowas}, there exists $\epsilon>0$ such that
\begin{equation}\label{2k}
K^{+}_\epsilon \cap K^{-}_\epsilon\subset V,
\end{equation}
where $K^{\pm}_\epsilon$ denote the $\epsilon$-neighborhood of
$K^{\pm}$. We define the compact set $D^{\pm} = \overline{D \setminus
K^{\pm}_\epsilon}$. Clearly, $D^{\pm}\cap K^{\pm}=\varnothing$, and
from (\ref{2k}) we have $D\subset D^+\cup D^-$.  From the uniform
convergence of $f^k(p)$ to $I^-$ on  $D^{-}$ we conclude that
there exists $N_2>0$ such that $f^{k}(p)\subset \hat V^-$ for
$k\ge N_2$ and $p\in D^-$. Analogously, $f^{-k}(p)\subset \hat
V^+$ for $p\in D^+$ and for $k\ge N_3>0$. Combining the above
considerations, we conclude that Claim 1 holds for $N=\max (N_1,
N_2, N_3)$.

\medskip

\noindent{\it Claim 2. Let $N$ be as in  Claim 1. Then $V_N^-\cup
V_N^+\cup V = \bC^n$.}

\medskip

We prove the claim by contradiction. We assume that there exists a
point $p \in \bC^n\setminus (V \cup V^{+}_N \cup V^{-}_N)$. By
Claim 1, either $f^{N}(p)\in \hat V^-$ or $f^{-N}(p)\in \hat V^+$.
If $f^{-N}(p)\in \hat V^+\cap \bC^n= V^{+}_0$, then $p\in
f^N(V^{+}_0)=V^{+}_N$ which is a contradiction. Now we suppose
that $f^N(p)\in \hat V^-$. If $f^N(p)\in V^{-}_0(N)$, then $p\in
f^{-N}(V^{-}_0(N))=V^{-}_N$, which is again a contradiction. The
remaining case is $f^{N} (p) \in (\hat V^- \cap \bC^n ) \setminus
V^{-}_{0} (N)$. Since $f^{N} (p)$ is not contained in $V^{-}_{0}
(N)$, it follows that for some $j\leq N$, $f^{N-j}(p)\in V^{+}_N$.
But then the invariance of $V^+_{N}$ under $f^{-1}$ implies that
$p\in V^+_{N}$. This proves the claim.

\medskip

We finally set $V^{\pm}=V_N^{\pm}$, where $N$ is such that $\bC^n
\setminus V \subset V^{+}_N \cup V^{-}_N $. We remark that $V^-$
cannot be empty, because otherwise there would exist a
neighborhood  of $I^-$ whose intersection with $\bC^n$ is
contained in $V^+$. Then, since the closure of $K^-$ contains $I^-$,
there exists a point $p\in K^-$ which is also contained in $V^+$.
By construction, $f^{-N}(p)\in \hat V^+$. But $\hat V^+\cap
K^-=\varnothing$ and this contradicts the fact that $K^-$ is an
$f$-invariant set. 

We redefine the set $V$ by setting $V=\bC^n\setminus(V^+\cup  V^-)$. 
Then $V$ is compact. Indeed, the union of $V^+$ and $V^-$ is not 
changed by the shrinking, and therefore, $V^+\cup V^-$ is open. 
We claim that $K\subset \inte V$. To show the claim, we note that by continuity for every $p\in K$ there exists  a neighborhood 
$U\subset \bC^n$ of $p$ such that $f^N(U)\cap \hat V^-=\emptyset$ and 
$f^{-N}(U)\cap \hat V^+=\emptyset$, which implies $U\subset V$.

It follows immediately from the construction that
$f^{-1}(V^+)\subset V^+$, i.e., inclusion
(iii) hold. To obtain property (i), we first observe that it is
sufficient to show $f(V^-_0(N))\subset V^-_0(N)$.
Let $p\in
V^{-}_0(N)$, in particular, $p\not\in f^N(V^+_N)$. Using the fact that
$f^{j}(V^+_N)\subset f^{j+1}(V^+_{N})$, we obtain
$f(p)\not\in f^{N+1}(V^+_N)\supset
f^{N}(V^+_N)$, which implies $f(p)\in V^-_0(N)$, and (i) holds.
We now show property (ii). If $p\in V^-$, then
$f(p)\in V^-$ by (i).  Let now $p\in V$, and assume $f(p)\in V^+$.
Then, it follows from (iii) that $p\in V^+$, which is a
contradiction. Analogously we obtain property (iv).
\end{proof}

\begin{corollary}\label{l1}
Let $f$ be a regular polynomial automorphism of $\bC^n$, and let
the compact set $V\subset \bC^n$ be defined as in Proposition
\ref{l1.5}. Then
\begin{equation}\label{eqwer}
\begin{array}{l}
f^{\pm 1}(K^\pm\cap V)\subset K^\pm\cap V,\\ f^{\pm 1}(J^\pm\cap
V)\subset J^\pm\cap V.
\end{array}
\end{equation}
\end{corollary}

\begin{proof}
Suppose $p\in K^+\cap V$. Then by Proposition \ref{l1.5}, $f(p)\in
V$ or $f(p)\in V^-$ . We only need to consider the case $f(p)\in
V^-$. It follows from the  construction that $f^{N+1}(p)\in
V^-_0(N)\subset \hat{V}^-$. On the other hand, $K^+$ is an
$f$-invariant set, and $K^+\cap \hat V^-=\varnothing$. This is a
contradiction. Similarly, one can easily verify the other cases.
The proof of the second inclusion in (\ref{eqwer}) follows from
the first inclusion and the $f$-invariance of $J^{\pm}$.
\end{proof}

\noindent
 {\it Remark. }We note that for $n=2$, that is, when $f$ is a
finite composition of generalized H\'enon mappings, $V$ can be
chosen to be  a closed bidisk of a sufficiently large radius R,
and $V^{-} = \{ (x,y)\in \bC^2 : |y|>R {\rm \ and \ } |y|>|x| \}
$, $V^{+} = \{ (x,y)\in \bC^2 : |x|>R {\rm \ and \ } |y|<|x| \} $,
(see \cite{BS}).

\medskip

For a set $X\subset \bC^n$ we define the stable and unstable sets
$W^s(X)$ and $W^u(X)$ as
\begin {equation}\label{equnstable}
\begin{array}{l}
W^s(X)=\{q\in \bC^n: {\rm dist}(f^k(q),f^k(X))\to 0\ {\rm as}\
k\to \infty\},\\ W^u(X)=\{q\in \bC^n: {\rm
dist}(f^{-k}(q),f^{-k}(X))\to 0\ {\rm as}\ k\to \infty\}.
\end{array}
\end{equation}

\begin{lemma}\label{l:bs}
Let $f$ be a regular polynomial automorphism of $\bC^n$. Then the
following hold.
\begin{enumerate}
\item[(i)] $W^s(K)=K^+$;
\item[(ii)] $W^u(K)=K^-$;
\item[(iii)] $\bigcup f^{k} (V^+) = \bC^n\setminus K^-$;
\item[(iv)] $\bigcup f^{-k} (V^-) = \bC^n\setminus K^+$;
\item[(v)] If $|\det Df|=1$, then $\inte K^+=\inte K^-=\inte K$;
\item[(vi)] If $|\det Df|<1$, then $K^-$ has zero Lebesgue measure,
  in particular, $\inte K^-=\emptyset$.
\item[(vii)] $\{f^{\pm k}: \ k\in\bN\}$ is a normal family on $\inte K^{\pm}$.
\end{enumerate}
\end{lemma}

\begin{proof}
(i) The inclusion  $W^s(K)\subset K^+$ follows from the
definition. The opposite inclusion follows from the fact that the
distance of $f^k(p)$ to $K$ converges uniformly to zero on compact
subsets of $K^+$. Applying the same arguments to $f^{-1}$ gives (ii).

(iii) We first assume $p\in\bC^n\setminus K^-$. Then $f^{-k}(p)$
converges to $I^+$, and therefore, $f^{-k}(p)\in V_0^+ \subset
V^+$ for sufficiently large $k$.  On the other hand, if $p\in
\bigcup f^{k} (V^+)$, then by Proposition \ref{l1.5} (i)
 there exists $k_0\in\bN$ such that $f^{-k}(p)\in V^+$ for all $k\geq k_0$.
Therefore, $f^{-k}(p)$ cannot converge to $K\subset \inte V$. Hence
$p\not\in K^-$, and property (iii) holds.

iv) Suppose $p\in\bC^n\setminus K^+$. Then $f^k(p)$ converge to
$I^-$ as $k\to\infty$. Let us show that $\bC^n\setminus K^+\subset
\bigcup f^{-k} (V^-)$. Suppose on the contrary that
$f^{k}(p)\not\in V^-$ for any $k>0$. Then, since the orbit of $p$
converges to $I^-$, there exists $k_0>0$ such that $f^k(p)\in V^+$
for all $k>k_0$. Let $N$ be as in Claim 1 of the proof of
Proposition \ref{l1.5}, and let $k> N$ be arbitrary. Then $f^k
(p)\in V^+$ implies $f^{k-N}(p)\in V^+_{0}$, i.e., all iterates of
$p$ stay in $V^+_0$. But $V^+_0\cap \hat{V}^-=\varnothing$ and
therefore, since $\hat{V}^-$ is a neighborhood of $I^-$, this
contradicts the fact that $f^k(p)$ converge to $I^-$. The opposite
inclusion can be proven similarly to case (iii).

(v) Assume on the contrary that there exists a ball
$B=B(p,r)\subset\inte K^+\setminus \inte K$. Without loss of
generality we may assume $B\subset \inte K^+\setminus K$, in
particular, $B\subset \bC^n\setminus K^-$. Since $f^{-k}$ converges
uniformly to $I^+$ on compact subsets of $\bC^n\setminus K^-$,
there exists a subsequence $(k_j)_{j\in\bN}$ such that the sets
$f^{-k_1}(B),f^{-k_2}(B),f^{-k_3}(B),\dots$ are pairwise disjoint.
Using that $f^{-1}$ is volume preserving we obtain that
$\vol(\inte K^+)=\infty$.  Here $\vol$ denotes the Lebesgue
measure in $\bC^n$. Thus there exists $r>0$ such that $\vol
(B(0,r)\cap \inte K^+)> \vol (V)$. By the uniform convergence of
$f^k$ on compact subsets of $ K^+$ there exists $k_0\in \bN$ such
that $f^{k_0}(\inte K^+\cap B(0,r))\subset V$. But this is a
contradiction to $\vol(f^{k_0}(\inte K^+\cap B(0,r)))=\vol(\inte
K^+\cap B(0,r))>\vol(V)$. Thus  $\inte K^+= \inte K$ holds. The
proof of the identity $\inte K^-= \inte K$ is analogous. 

Property (vi) follows analogously to the case $n=2$ (see \cite{FM}).
Finally, (vii) follows from \cite[Prop. 2.2.7]{Si}.
\end{proof}

%%_________________________ section_____________________________

\section{Hyperbolicity}

For generalized H\' enon maps in $\bC^2$  the concept of hyperbolicity
was studied in detail by Bedford and Smillie (see \cite{BS}).
Using the filtration properties obtained in Section 3, we
generalize in this section  some of the results of \cite{BS} to
regular polynomial automorphisms of $\bC^n$.

We first give some basic definitions. We refer to \cite{KH} for
details. Let $f$ be a  regular polynomial automorphism  of
$\bC^n$. We say that a compact $f$-invariant set $\Lambda\subset
\bC^n$ is a {\it hyperbolic set} for $f$  if there exists a continuous
$Df$-invariant splitting of the tangent bundle $T_\Lambda \bC^n=
E^u\oplus E^s$ such that $Df$ is uniformly expanding on $E^u$ and
uniformly contracting on $E^s$.

An important feature of hyperbolic sets is that we can associate
with each point $p\in \Lambda$ its local unstable/stable manifold
$W^{u/s}_\epsilon(p)$.  The local unstable/stable manifolds are
complex manifolds of the same (complex) dimension as $E_p^{u/s}$.
We denote by $W^{u/s}(p)$ (global) unstable/stable manifold of $p$
(see \eqref{equnstable}). It follows from the work of Jonsson and
Varolin \cite{JV} that for all $p$ in a set   of total probability
the global unstable/stable manifolds $W^{u/s}(p)$ are
biholomorphic copies of  $\bC^k$ where $k=\dim_\bC E^{u/s}_p$. We
call $\dim_\bC E^{u/s}_p$ the unstable/stable index of $\Lambda$
at $p$. Note that unstable/stable index is locally constant. If
$\Lambda$ is a hyperbolic set for $f$, we say that $\Lambda$ is
{\it locally maximal} if there exists a neighborhood $U$ of $\Lambda$
such that every hyperbolic set of $f$ in $U$ is contained in $\Lambda$.
We say that an $f$-invariant set $X$ has a {\it local product structure}
if for all $p,q\in X$ we have $W^s(p)\cap W^u(q)\subset X$. We now
consider the situation when $J$ is a hyperbolic set for $f$.

\begin{proposition}\label{l:hyper}
Let $f$ be a regular polynomial automorphism of $\bC^n$, and suppose
that $J$ is a hyperbolic set for $f$. Then the following holds.
\begin{enumerate}
\item[(i)] If $p\in J$, then $W^{s/u} (p)\subset J^\pm$;
\item[(ii)] The set $J$ has a local product structure;
\item[(iii)] The set $J$ is locally maximal, and $W^{u/s}(J)=\bigcup_{p\in J}W^{u/s}(p)$;
\item[(iv)] $W^{s/u}(J)\subset J^\pm$.
\end{enumerate}
\end{proposition}

\begin{proof}
(i). Without loss of generality we show the inclusion only for
$W^s(p)$ The proof for $W^u(p)$ is analogous. Clearly,
$W^s(p)\subset K^+$. Suppose there exists a point $q\in W^s(p)\cap
{\inte } K^+$. Then, since the family $\{f^k: k\in\bN\}$ is normal
in a neighborhood of $q$, the derivatives of $f^k$ at $q$ are
bounded. On the other hand, by extending the hyperbolic structure
of $f$ to a neighborhood of $J$, it follows that the derivatives
of $f^k$ at $q$ are unbounded.

(ii) If $p,q\in J$, then by (i), $W^s (p)\in J^+$ and $W^u(q)\in J^-$,
therefore, the intersection is in $J$.

(iii) The local product structure combined with hyperbolicity  implies that $J$ is  locally maximal
(see \cite[Prop. 8.22]{Sh}). The second statement is an application of the
shadowing lemma for locally maximal hyperbolic sets (see \cite{Bo}).

Finally, (iv)  follows from (i) and (iii).
\end{proof}

Let $C$ be a connected component  of $\inte K^+$. We say that $C$
is  {\it periodic}, if there exists $N\in \bN$ such that $f^N(C)=C$.
Otherwise we call $C$  {\it wandering}. If
 $\alpha$ is a periodic point such that for all $p$ in a neighborhood of $\alpha$   we have
$f^k(p)\to f^k(\alpha)$ as $k\to\infty$, then we call $\alpha$ an
{\it attracting periodic point}. Furthermore, $C=\{p\in \bC^n: f^k(p)\to
f^k(\alpha)\}$ is a periodic connected component of $\inte K^+$
and is called  the {\it basin of attraction} of  $\alpha$.

\begin{theorem}\label{k+}
Let $f$ be a regular polynomial automorphism of $\bC^n$ with $|\det
Df|\leq 1$. Suppose
$J$ is a hyperbolic set for $f$. Then  the following holds.
\begin{enumerate}
\item[(i)] There are no wandering components in $\inte K^+$;
\item[(ii)] Each periodic component of $\inte K^+$ is the basin of
  attraction of an attracting
  periodic point;
\item[(iii)] There are at most finitely many basins of attraction.
\end{enumerate}
\end{theorem}

The proof of Theorem \ref{k+} is analogous to that of Theorem~5.6 in
\cite{BS}. We note that the references to Propositions~5.1 and 5.2 in
the proof of \cite{BS} should be replaced by references to
Proposition~\ref{l:hyper} stated above.

\begin{corollary}
Let $f$ be a regular polynomial automorphism of $\bC^n$. Assume $J$ is
hyperbolic and
$|\det Df|=1$. Then $\inte K^+ = \inte K^- =\inte  K =
\varnothing $.
\end{corollary}

\begin{proof}
By Theorem~\ref{k+}, the interior of $K^+$ is a finite union of
basins of attraction, but since $|\det Df|=1$, it is impossible to
have a basin of attraction. Hence $\inte K^+ = \varnothing$.
Therefore, the corollary follows from Lemma~\ref{l:bs} (v).
\end{proof}

We need the following definitions. Let $X$ be a topological space
and let $T:X\to X$ be a continuous map. We call $x\in X$ a {\it
nonwandering} point of $T$ if for every neighborhood $U$ of $x$
there exists $k\in\bN$ such that $U\cap T^k(U)\not=\emptyset$.
Otherwise we call $x$ {\it wandering}. The  set of nonwandering points
of $T$ is called the {\it nonwandering set} of $T$ and is denoted
by $\Omega(T)$.

We say that a regular polynomial automorphism $f$ of $\bC^n$ is
{\it hyperbolic} if $J$ is a hyperbolic set for $f$ and the
periodic points of $f|_J$ are dense in $J$. We note that this
definition of hyperbolicity is equivalent to $J$ being hyperbolic
and $\Omega(f|J)=J$ (see for instance \cite{KH}).

\begin{theorem}\label{thghj}
Let $f$ be a hyperbolic regular polynomial automorphism  of
$\bC^n$ with $|\det Df|\leq 1$.  Then
\begin{enumerate}
\item[(i)] $W^s(J)=J^+$;
\item[(ii)] $W^u(J)=J^-\setminus  \{\alpha_1,\dots,\alpha_m\}$, where
  the $\alpha_i$ are the attracting periodic points of $f$;
\item[(iii)] $\Omega(f) = J\cup  \{\alpha_1,\dots,\alpha_m\}$.
\end{enumerate}
\end{theorem}

\begin{proof}
(i) To prove $W^s(J)=J^+$ we observe that by
Proposition~\ref{l:hyper} (iv) we have $W^s(J)\subset J^+$. To
show the reverse inclusion we notice that it follows from
Lemma~\ref{l:bs} (i) that $J^+\subset W^s(K)$. If $p\in J^+$, then
the iterates of $p$ converge to $K\cap J^+=K^-\cap J^+$. To prove
(i) we claim that $K^-\cap J^+=J^-\cap J^+=J$. In the case  $|\det
Df|=1$ the claim follows from Lemma \ref{l:bs} (v).  For $|\det
Df|<1$ the claim follows from  Lemma \ref{l:bs} (vi).

(ii)  Obviously every attracting period point belongs to $\inte
K^+$, and therefore, $W^u(J)\subset J^-\setminus
\{\alpha_1,\dots,\alpha_m\}$ follows from
Proposition~\ref{l:hyper} (iv). In order to show the opposite
inclusion we consider $p\in J^-\setminus
\{\alpha_1,\dots,\alpha_m\}$. If $p\in \bC^n\setminus K^+$, then
the backward orbit of $p$ converges to $K\cap \partial
(\bC^n\setminus K^+)=K\cap J^+$. Using the fact that $J^-$ is a
closed invariant set, we conclude that the backward orbit of $p$
must converge to $J^+\cap J^-=J$. If $p\in J^+$, then $p\in  J$,
and there is nothing to prove. To complete the proof of (ii) we
have to consider the case  $p\in \inte K^+$. By Theorem \ref{k+}, 
there exists an attracting periodic point $\alpha_i$ such
that $p$ is contained in the basin of attraction $C$ of $\alpha_i$. 
Without loss of generality we assume that $\alpha_i$ is an attracting 
fixed point, i.e., $C$ is an $f$-invariant
component. Let $V$ be as in Proposition \ref{l1.5} and let $U\subset
V\cap C$ be an open neighborhood of $\alpha_i$ such that
$f(\overline{U})\subset U$. Such a  set $U$ always exists since
$\alpha_i$ is an attracting fixed point. Obviously $\bigcup
f^{-k}(U)$ is an exhaustion of $C$. Thus $\bigcup
f^{-k}(U\cap J^-)$ is an exhaustion of $C\cap J^-$. Therefore, since $p\in K$, we may conclude  that
the backward orbit of $p$ cannot have a cluster point in $C$, and thus
it must converge to $\partial C\subset J^+$. Thus $p\in J^-\setminus
\{\alpha_1,\dots,\alpha_m\}$ implies $p\in W^u(J)$.

(iii) Evidently every periodic point of $f$ belongs to the nonwandering
set of $f$. Since the nonwandering set is closed, it follows that
 $\Omega(f)\supset J\cup  \{\alpha_1,\dots,\alpha_m\}$.
  Let  $p\in \bC^n$ be a nonwandering
point for $f$. If $p$ is not an attracting periodic orbit, then $p$
can not  belong to $\inte K^+$,  because otherwise Theorem \ref{k+} would imply
that its forward orbit
converges to an attracting periodic orbit. On the other hand
$p\not\in \bC^n\setminus K^+$, because in this case the forward orbit
of $p$ would
converge to $I^-$. Hence
$p\in J^+$. It follows from (i) that $p\not\in J^+\setminus J$.
This implies $p\in J$, which completes the proof.
\end{proof}

We say that a diffeomorphism on a Riemannian manifold  is 
Axiom~A  if its nonwandering set is a
hyperbolic set and the periodic points  are dense in the
nonwandering set.

\begin{corollary}\label{coraxioma}
Let $f$ be a hyperbolic regular polynomial automorphism of
$\bC^n$. Then $f$ is Axiom A.
\end{corollary}
\begin{proof}
This is a consequence of Theorem \ref{thghj}~(iii) and the fact
that every attracting periodic point is an isolated hyperbolic
set.
\end{proof} \noindent {\it Remarks.} 
\begin{enumerate}
\item
As it was noted in the
introduction, if $n=2$ and $J$ is hyperbolic, then it follows that
the periodic points are dense in $J$, see \cite{BS}. This provides
a weaker definition of hyperbolic maps. We do not know whether the
analogous result holds in the case $n>2$. However, there exist
examples of diffeomorphisms of higher dimensional real manifolds
with the property that the non-wandering set is a hyperbolic set,
but the periodic points are not dense in it, see \cite{D}.
\item
It is shown in \cite[Thm. 9.3.14]{MNTU} that in the case $n=2$ 
hyperbolicity is equivalent to Axiom~A. We do not know whether the 
analogous result holds for $n>2$. The difficulty is to prove that the Julia set of an Axiom~A regular polynomial map 
is a subset of the nonwandering set. This is shown in the case $n=2$
with methods that are not available in higher dimensions.

\end{enumerate}
%%--------------------------- dimension stuff  ------------------------

\section{Dimension theory}
 In this section we study  the
Hausdorff dimension and box dimension of the Julia sets of a regular polynomial automorphism of $\bC^n$.

\subsection{The general case}
We first consider the dimensions of the Julia sets without the assumption 
of hyperbolicity.

Let $f$ be a regular polynomial automorphism of $\bC^n$ and let
$V\subset \bC^n$ be a compact set with $K\subset \inte V$ and $f^{\pm
  1}(J^\pm\cap V)\subset J^\pm\cap V$ (see Corollary \ref{l1}). We define
\begin{equation}\label{defspm}
 s^\pm_V=\lim\limits_{k\to\infty}\frac{1}{k}\log\left(\max\{\lVert Df^{\pm k}(p)\rVert: p\in
J^\pm\cap V\}\right).
\end{equation} The submultiplicativity of the
operator norm guarantees the existence of the limit defining
$s^\pm_V$.
Since all norms in $\bC^n$ are equivalent, the value of
$s^\pm_V$ is independent of the norm. Moreover, since the saddle points of $f$ are contained in $J\subset J^\pm\cap
V$ (see \cite{Si}), the quantity $s^\pm$ is strictly positive.

\begin{lemma}\label{l2}
The value of $s^\pm_V$ is  independent of the choice of $V$.
\end{lemma}

\begin{proof} It is shown in Lemma \ref{l:bs} that
$W^s(K)=K^+$ and $W^u(K)=K^-$. Therefore, the proof
follows by a standard argument.\end{proof}

In view of Lemma \ref{l2} we set $s^\pm=s^\pm_V$. Given a set $A\subset
\bC^{n}\cong \bR^{2n}$ we denote by $\dim_H A$ the Hausdorff dimension of
$A$ and (provided $A$ is bounded) by $\overline{\dim}_B A$
its upper box-dimension   (see \cite{M} for details). Then
$\dim_H A \leq \overline{\dim}_B A$ holds for an arbitrary set $A$,
while the equality holds if $A$ is a sufficiently regular set. We now consider the volume decreasing case, i.e., when $|\det Df|<1$. The
following theorem provides an upper bound for the dimension of $K^-$.

\begin{theorem}\label{th4-}
 Let $f$ be a volume-decreasing regular polynomial
 automorphism of $\bC^n$. Assume $V$ is as in Corollary \ref{l1}.
Then
\begin{equation}\label{eq4-}
\overline{dim}_B K^-\cap V \leq 2n +
\frac{2\log|\det Df|}{s^-}<2n.
\end{equation}
\end{theorem}

\begin{proof} By Lemma \ref{l:bs} (vi),  $K^-=J^-$. Therefore, it is
  sufficient to show inequality
\eqref{eq4-} for $J^-\cap V$. Note that the real Jacobian of $f^{-1}$ as
a map of $\bR^{2n}\cong \bC^n$ is equal to $|\det Df|^{-2}$.
The result follows now immediately from \cite[Thm. 1.1]
{W1}.\end{proof}

\noindent{\it Remark.}
Since $W^u(K)=K^-$ we can
define an exhaustion $V_k=f^k(V\cap K^-)$ of $K^-$. This implies
that the upper bound in inequality \RE{eq4-} provides also an upper
bound for the Hausdorff dimension of $K^-$.

\begin{corollary}\label{corkk4}
Let  $f$ be regular polynomial automorphism of $\bC^n$ which is not
volume-preserving. Then $\overline{\dim}_B K <2n$.
\end{corollary}

\begin{proof}
If $f$ is not volume-preserving, then either $f$ or $f^{-1}$ is
volume-decreasing. The result follows immediately from Theorem \ref{th4-}.
\end{proof}
\noindent {\it Remark. }
It should be noted that Corollary \ref{corkk4} does not hold
without the assumption that $f$ is not volume-preserving. In fact,
there exist volume-preserving regular polynomial automorphisms
having a Siegel ball, in which case $K$ has a non-empty interior.

\medskip

For a continuous map $T$ on a compact metric space $X$ we denote
by $\htop(T)$ the topological entropy of $T$ (see \cite{Wa} for
details). We now show that the Julia set $J$ carries the full
entropy of $f$.
\begin{theorem}\label{prohtop}
Let  $f$ be regular polynomial automorphism of $\bC^n$. Assume
that $I^-$ has dimension $l-1$. Then $h_{\rm top}(f|J)=l\log d$.
\end{theorem}
\begin{proof}
Without loss of generality we assume $|\det Df|\leq 1$. We claim
that $K\setminus J\subset \inte K^+$. Indeed, we have $K=(\inte K^+\cup
J^+)\cap (\inte K^-\cup J^-)$. This implies

\begin{equation}\label{eqhoff}
K\setminus J\subset \inte K^+\cup \inte K^-.
\end{equation}
 If $|\det Df|<1$, then by Lemma \ref{l:bs} (vi), $\inte
 K^-=\emptyset$. Therefore,  \eqref{eqhoff}
implies $K\setminus J\subset \inte K^+$. On the other hand, if
$|\det Df|=1$, then by Lemma \ref{l:bs} (v), $\inte K^+= \inte
K^-= \inte K$. Again by \eqref{eqhoff} we obtain
$K\setminus J\subset \inte K^+$.

Sibony observed in \cite{Si} that $\htop(f|K)= l\log d$. Let
$\epsilon>0$. It follows from the variational principle (see for
instance \cite{Wa}) that there exists an $f$-invariant probability
measure $\mu$ on $K$ with $h_\mu(f|K)> l\log d-\epsilon$, where
$h_\mu(f|K)$ denotes the measure-theoretic entropy of $f|K$ with
respect to $\mu$. Let $\tau$ be an ergodic decomposition of $\mu$.
This means that $\tau$ is a probability measure on the metrizable
space $\cM$ of $f$-invariant probability measures on $K$ which
puts full measure on the subset $\cM_E$ of ergodic measures.
Furthermore,
\begin{equation}\label{ergdec}
\int_{\cM}\int_K \varphi\, d\nu\, d\tau(\nu) = \int_K\varphi\,
d\mu
\end{equation}
for every $\varphi\in C(K,\bR)$. Since
\begin{equation}
h_\mu(f|K)=\int_{\cM} h_\nu(f|K)\,d\tau(\nu),
\end{equation}
there exists $\nu\in \cM_E$ such that $h_\nu(f|K)> l\log
d-\epsilon$.

Next we  claim that ${\rm supp\ } \nu\subset J$. If not, then
$\nu(K\setminus J)>0$. Since $\nu$ is ergodic, this would imply
${\rm supp\ }\nu\subset K\setminus J\subset\inte K^+$. By Ruelle's inequality $\nu$ must have at least one positive Lyapunov exponent.
On the other hand, since $\nu$ has compact support, Lemma \ref{l:bs}
(vii) implies that the derivatives of $f^k$ are bounded on ${\rm
supp}\ \nu$. This contradicts the existence of a positive Lyapunov
exponent. Hence ${\rm supp\ } \nu\subset J$.

Finally, since $\epsilon$ can be chosen arbitrarily small, the
variational principle yields $\htop (f|J)=l\log d$.
\end{proof}
The following is an immediate consequence of the proof of Theorem
\ref{prohtop}.
\begin{corollary}\label{corhkj}
Let  $f$ be a regular polynomial automorphism of $\bC^n$. Then
$h_{\rm top}(f|K\setminus J)=0$.
\end{corollary}
We note that in Corollary \ref{corhkj}, $\htop$ denotes the
entropy for maps on non-compact spaces (see \cite{Wa}).

Another consequence of Theorem \ref{prohtop} is that the upper box
dimension of $J$ is strictly  positive.
\begin{corollary}\label{corkk3}
Let $f$ be a regular polynomial automorphism of $\bC^n$. Assume
that $I^-$ has dimension $l-1$. Define
\[
 s^\pm=\lim\limits_{k\to\infty}\frac{1}{k}\log\left(\max\{\lVert Df^{\pm k}(p)\rVert : p\in
J\}\right).
\]
Then
\begin{equation}\label{eqbospos}
\overline{\dim}_B J \geq \max\left\{\frac{l\log
d}{s^+},\frac{l\log
      d}{s^-}\right\}.
\end{equation}
In particular, $\overline{\dim}_B J>0$.
\end{corollary}
\begin{proof} We have $\htop (f|J)=\htop(f^{-1}|J)=l\log
  d$. Therefore, inequality
\eqref{eqbospos} follows from \cite[Thm. 3.2.9]{KH} and a standard
limit argument. \end{proof}

Next we give a lower bound for the Hausdorff dimension of $J^+$ in
the case when $I^-$ has  dimension zero. For this we introduce the
positive Green function
\begin{equation}
G^{+}(p) = \lim_{k\to\infty} \frac{1}{d^k} \log^+|f^{ k}(p)|.
\end{equation}
We note that $G^{+}$ is a well-defined   H\"older continuous
function with $K^{+}=\{G^{+}=0\}$ (see \cite{Si} for details) .

\begin{proposition}\label{prosibony}
Let $f$ be a regular polynomial automorphism of $\bC^n$. Assume that
$I^-$ has dimension $0$ and let $s^+$ be as in \eqref{defspm}.
Then for all $s^+_0
> s^+$ the positive Green function $G^+$ is H\"{o}lder continuous on  compact subsets of $\bC^n$ with H\"{o}lder exponent $\frac{\log
d}{s^+_0}$. Furthermore,
\begin{equation} \dim_H J^+\geq 2n-2+\frac{\log d}{s^+}>2n-2.\label{ghjk}
\end{equation}
\end{proposition}
\begin{proof}
Let $s^+_0>s^+$. Using the filtration properties (see Proposition
\ref{l1.5}), one can show the H\"older continuity of $G^+$ on
compact subsets  of $\bC^n$ with H\"older exponent $\frac{\log
d}{s^+_0}$ analogously as it was done by Forn\ae ss and Sibony
\cite{FS} in the case $n=2$ (see also \cite{Si}). By \cite[Prop.
2.2.10]{Si} the positive Green function $G^+$ is pluriharmonic on
$\bC^n\setminus J^+$. Furthermore, the maximum principle for
pluriharmonic functions implies that $G^+$ cannot be extended as a
pluriharmonic map to any neighborhood of any point of $J^+$.
Inequality \eqref{ghjk} follows now by a classical result of
Carleson about the Hausdorff dimension of removable sets for
H\"older continuous harmonic functions (see \cite{C}).
\end{proof}
\noindent
{\it Remarks.}
\begin{enumerate}
\item[(i)] The analogous result holds for the Hausdorff dimension
of $J^-$ if $I^+$ is zero-dimensional.
\item[(ii)]
We note that Proposition \ref{prosibony} is only of interest if
$\inte K^+=\emptyset$, because otherwise $J^+$ has topological
dimension $2n-1$.
\end{enumerate}

\subsection{Hyperbolic maps}
We now consider  hyperbolic maps. It is well-known that a locally
maximal hyperbolic set, which carries positive topological
entropy, has positive Hausdorff dimension. For a hyperbolic
regular automorphism $f$ of $\bC^n$ the positivity of the
Hausdorff dimension of the Julia set can be shown for instance by
the following argument: Let $J=J_1\cup\dots\cup J_m$ be the
decomposition of $J$ into basic sets (see e.g. \cite{Bo} for
details); We note that in the case $n=2$ the Julia set $J$ is the unique
basic set of $f$ which is not an attracting periodic orbit (see \cite{BS}). Since $\htop (f|J)=l\log d$ (see
Theorem \ref{prohtop}), there exists $i\in\{1,\dots,m\}$ such that
$\htop(f|J_i)=l\log d$. Here $l=\dim I^- + 1$. It is a well-known
fact that there exists a unique $f$-invariant probability measure
$\mu_i$ of maximal entropy for $f|J_i$ (see for example
\cite{KH}). Moreover $\mu_i$ is ergodic. We define
\begin{equation}\label{dimJpos}
s_i^\pm=\lim_{k\to\infty} \frac{1}{k}\log\left( \max\{\lVert
Df^{\pm k}(p) \rVert: p\in J_i\}\right).
\end{equation}
Applying \cite[Cor. 5.1]{Y} yields
\begin{equation}
 l\log d\left(\frac{1}{s^+_i}+\frac{1}{s^-_i}\right)\leq \dim_H \mu_i \leq
\dim_H J_i\leq \dim_H J,
\end{equation}
where $\dim_H \mu_i=\inf\{\dim_H A: \mu_i(A)=1\}$ denotes the
Hausdorff dimension of the measure $\mu_i$. We note that in
general, the Hausdorff dimension of an $f$-invariant measure
provides only a rough estimate of the Hausdorff dimension of the
Julia set. In fact, it was shown in \cite{W3} that for a generic
hyperbolic polynomial automorphism $f$ of $\bC^2$ there exists
$\epsilon>0$ (which depends on $f$) such that $\dim_H \nu < \dim_H
J -\epsilon$ for all ergodic $f$-invariant probability measures
$\nu$.

Let $f$ be a  hyperbolic regular polynomial automorphism of
$\bC^n$, let $J_i\subset J$ be a basic set of $f$ and $\varphi\in
C(J_i,\bR)$. We denote by $P(f|J_i,\varphi)$ the topological
pressure of $\varphi$ with respect to $f|J_i$ (see \cite{KH} for
the definition and details). We consider the function $\phi^{u}=-
\log \lVert Df|E^{u}\rVert$. Note that $\phi^{u}$ is H\"older
continuous, see \cite{Bo}.

We now consider the case when the unstable index $J_i$ is
identically one.
\begin{theorem}\label{thboru}
Let $f:\bC^n\to\bC^n$ be a hyperbolic regular polynomial
automorphism and let $J_i\subset J$ be a basic set of $f$. Assume that
the unstable index of $J_i$ is identically one. Then $t^u=dim_H
W^u_\epsilon(p)\cap J_i$ is independent of $p\in J_i$ and
$0<t^u<2$. Moreover, $t^u$ is given by the unique solution of
\begin{equation}\label{eqbowen}
P(f|J_i,t\phi^u)=0.
\end{equation}
\end{theorem}

Equation \eqref{eqbowen} is usually called the Bowen--Ruelle
formula. We refer to $t^u$ as the Hausdorff dimension of the
unstable slice. Theorem \ref{thboru} is a special case of
\cite[Thm. 22.1]{P}. In the case $n=2$, Theorem \ref{thboru} is
due to Verjovsky and Wu \cite{VW}. We note that in this situation
the stable index of the basic set $J$ is also identically one, and
we obtain the analogous result to Theorem \ref{thboru} for the
Hausdorff dimension of the stable slice $t^s=\dim_H
W^s_\epsilon(p)\cap J.$ Moreover, $\dim_H J = t^u+t^s$ (see
\cite{W2} and the references therein).

The following result is a version of \cite[Thm. 4.1]{W2}. The proof is
analogous.

\begin{theorem}\label{thint}
Let $f$ be a hyperbolic regular polynomial automorphism of $\bC^n$
and let $J_i\subset J$ be a basic set of $f$  . Assume that the
unstable index of $J_i$ is identically one. Then $\dim_H
W^s(J_i)=t^u+2n-2$. In particular, $2n-2<\dim_H W^s(J_i)<2n$.
\end{theorem}

Let now $A\subset \bC^k$ 
be an open set and let $\{f_a: a\in A\}$
be a holomorphic family of hyperbolic regular polynomial
automorphisms of $\bC^n$. Let $a_0\in A$ and let $J_{a_0,i}\subset
J_{a_0}$ be  a basic set of $f_{a_0}$. Let $U\subset \bC^n$ be a
neighborhood of $J_{a_0,i}$ with the property that for all $a\in
A$ close enough to $a_0$, $f_a$ has a basic set $J_{a,i}\subset U$
such that $f_{a_0}|J_{a_0,i}$ is conjugate to $f_a|J_{a,i}$. For $p\in
J_{a,i}$ we
denote by $t^u_a$ the Hausdorff dimension of $W^u_\epsilon(p)\cap
J_{a,i}$. Recall that by Theorem \ref{thboru} $t^u_a$ does not depend
on $p$. The following result can be proven analogously to the
corresponding results in the case $n=2$ (see \cite{VW}, \cite{W2}).

\begin{theorem}\label{thhou}
Assume that the unstable index of $J_{a,i}$ is identically one  in a
neighborhood of $a_0\in A$.  Then the functions $a\mapsto t_a^u$
and $a\mapsto \dim_H W^s(J_{a,i})$ are real-analytic and
plurisubharmonic in a neighborhood of $a_0\in A$.
\end{theorem}
\noindent
 {\it Remark. } The corresponding versions  of Theorems \ref{thboru},
 \ref{thint}, \ref{thhou} for the stable slices also hold provided
 the stable index has dimension identically one.

Let $f$ be a hyperbolic regular polynomial automorphism of $\bC^n$
and let $J_i\subset J$ be a basic set of $f$. We define
$\varphi^{u/s}:J_i\to \bR$ by $\varphi^{u/s}(p)=\mp\log
\vert\lambda(p)\vert$, where $\lambda(p)$ denotes the Jacobian of
the linear map $Df^{\pm 1}(p)|E^{u/s}_p$. The following is the
main result of this section.

\begin{theorem}\label{thmaindim}
Let $f$ be a hyperbolic regular polynomial automorphism of $\bC^n$
and let $J_i\subset J$ be a basic set of $f$. Define
\begin{equation*}
W^{s/u}_\epsilon (J_i)=\bigcup_{p\in J_i} W^{s/u}_\epsilon(p),
\end{equation*}
\begin{equation*}
s^{\pm}=\lim_{k\to\infty} \frac{1}{k}\log\left( \max\{\lVert
Df^{\pm k}(p) \rVert: p\in J_i\}\right).
\end{equation*}
Then
\begin{equation}\label{equpbound1}
\overline{\dim}_B W^{s/u}_\epsilon (J_i)  \leq 2n +
\frac{P(f|J_i,\varphi^{u/s})}{s^{\pm}}< 2n.
\end{equation}
\end{theorem}

\begin{proof}
We prove the result only for $W^{s}_\epsilon (J_i)$. The proof for
$W^{u}_\epsilon (J_i)$ is entirely analogous. Since $J_i\subset
J$, its unstable index must be at least one,
 which implies $s^+>0$. By Proposition \ref{l:hyper},
 $W^s_\epsilon(J_i)\subset J^+$, in particular, $W^s_\epsilon(J_i)$ is
 not a neighborhood of $J_i$.
Therefore, we may conclude from  \cite[Prop. 3.10, 4.8, Thm. 4.11]{Bo}  that
 $P(f|J_i,\varphi^u)<0$. This gives the inequality on the right.

Let $\delta>0$. It follows by a simple continuity argument that there
exist $\epsilon>0$ and $k_\delta\in \bN$ such that for all
$p\in B(W^{s}_\epsilon (J_i),\epsilon)=\{p\in \bC^n:\exists
q\in W^{s}_\epsilon (J_i),\  |p-q|<\epsilon\}$ we have
\begin{equation}
\lVert Df^{k_\delta}(p)\rVert < \exp(k_\delta(s^++\delta)).
\end{equation}
From now on we consider the map $g=f^{k_\delta}$. Note that
$J_i$ is also a basic set of $g$. Evidently $W^{s}_\epsilon (J_i)$
is forward invariant under $g$. It follows from the variational principle
that $P(g|J_i,\varphi^u)=k_\delta P(f|J_i,\varphi^u)$; moreover
$s^+_g=k_\delta s^+_f$. It is thus sufficient to prove the
left-hand side of inequality \eqref{equpbound1} for $g$. Let $p\in
J_i$ and $k\in\bN$. Following \cite{Bo}, we
define
\begin{equation}\label{eqbo1}
B(p,\epsilon,k)=\{q\in \bC^n: |g^i(p)-g^i(q)|<\epsilon, i=0,\dots,k-1\},
\end{equation}
and $B(J_i,\epsilon,k)=\bigcup_{p\in J_i} B(p,\epsilon,k)$. Making
$\epsilon$ smaller if necessary, it follows from \cite[Prop. 4.8]{Bo}
that
\begin{equation}\label{eqbo2}
P(g|J_i,\varphi^u)=\lim_{k\to\infty}\frac{1}{k}\log(
\vol(B(J_i,2\epsilon,k))).
\end{equation}
For simplicity we write
$b=P(g|J_i,\varphi^u)$. From \eqref{eqbo2} we obtain that if $k$
is sufficiently large, then
\begin{equation}\label{eqbo4}
\vol(B(J_i,2\epsilon,k))< \exp(k(b+\delta)).
\end{equation}
For all $k\in \bN$ we define real numbers
\begin{equation}
r_k = \frac{\epsilon}{\exp(s^++\delta)^k}
\end{equation}
and neighborhoods $B_k = B(W^{s}_\epsilon (J_i),r_k)$ of
$W^{s}_\epsilon (J_i)$. Let $q\in B_k$. Then there exists $ p\in
W^{s}_\epsilon (J_i)$ with $|p-q| < r_k$. An elementary induction
argument in combination with the mean-value theorem implies
$|g^{i}(p)- g^{i}(q)| < \epsilon$ for all $i\in\{0,\dots,k-1\}$.
Since $p$ is contained in the local stable manifold of size $\epsilon$
of a point in
$J_i$  it follows that $q\in
B(J_i,2\epsilon,k)$. Hence $B_k\subset B(J_i,2\epsilon,k)$.
Therefore, \eqref{eqbo4} implies
\begin{equation}
\vol(B_k)<\exp(k(b+\delta))
\end{equation}
for sufficiently large $k$.
Let us recall that for $t\in [0,2n]$ the  $t$-dimensional upper
Minkowski content of a bounded set $A\subset \bC^n\cong\bR^{2n}$ is defined by
\begin{equation}
M^{*t}(A)= \limsup_{\rho\to 0} \frac{\vol(A_\rho)}{(2\rho)^{2n-t}},
\end{equation}
where $A_\rho=\{p\in\bC^n:\exists q\in A: |p-q|\leq \rho\}$. Let
$t\in [0,2n]$ and $\rho_k={\textstyle\frac{r_k}{2}}$ for all
$k\in\bN.$ Then we have
\begin{equation}\label{glej-4}
\begin{array}{lcl}
M^{*t} (W^s_\epsilon(J_i))&=&  \limsup\limits_{\rho\to 0}
\frac{\textstyle
\vol(W^s_\epsilon(J_i)_\rho)}{\textstyle(2\rho)^{2n-t}}\\ &\leq&
\limsup\limits_{k\to\infty } \frac{\textstyle
\vol(W^s_\epsilon(J_i)_{\rho_k})}{\textstyle(2\rho_{k+1})^{2n-t}}\\
&\leq &
 \limsup\limits_{k\to \infty} \frac{\textstyle \vol(B_k)}{\textstyle (r_{k+1})^{2n-t}}\\
&\leq & \frac{\textstyle
\exp(s^++\delta)^{2n-t}}{\textstyle
\epsilon^{2n-t}}\lim\limits_{k\to\infty}\left(\textstyle\exp(
s^++\delta)^{2n-t}\textstyle \exp(b+\delta)\right)^k.
\end{array}
\end{equation}
Let  $t > 2n + \frac{b+\delta}{ s^++\delta}$.  Then
$\exp(s^++\delta)^{2n-t} \exp(b+\delta) < 1.$ This implies
$M^{*t}(W^s_\epsilon(J_i)) = 0$, in particular, $t\geq
\overline{\dim}_B W^s_\epsilon(J_i)$ (see \cite{M}). Since $\delta$
can be chosen arbitrarily small, the result follows.
\end{proof}

\begin{corollary}\label{corjk2n}
Let $f:\bC^n\to\bC^n$ be a hyperbolic regular polynomial automorphism.
Then $\overline{\dim}_B J<2n$.
\end{corollary}
\begin{proof}
By the spectral decomposition, $J$ is the union of  finitely many basic
sets. By Theorem \ref{thmaindim}, each of these basic sets has upper box
dimension strictly smaller than $2n$, and the result follows.
\end{proof}
\noindent {\it Remark. } We note that Corollary \ref{corjk2n} is
only of interest if $f$ is volume-preserving, because otherwise
$\overline{\dim}_B J < 2n$ holds even without the assumption of
hyperbolicity (see Corollary \ref{corkk4}).
\begin{theorem}\label{jpmk4}
Let $f$ be a hyperbolic regular polynomial automorphism of $\bC^n$
and let $J_1,...,J_m$ be the basic sets of $f$ which are contained
in $J$. For each $i\in\{1,...,m\}$ define
\[
s^\pm_i=\lim_{k\to\infty} \frac{1}{k}\log\left( \max\{\lVert
Df^{\pm k}(p) \rVert: p\in J_i\}\right)
\]
and $b_i^\pm=P(f|J_i,\varphi^{u/s})$. Then
\begin{equation}\label{eqghj}
\dim_H J^\pm\leq 2n+\max\{b^\pm_i/s^{\pm}_i\}<2n.
\end{equation}
\end{theorem}
\begin{proof}
Without loss of generality we show the result only for $J^+$. By
Theorem \ref{thmaindim} we have
\begin{equation}
\overline{\dim}_B \left(\bigcup_{p\in J_i}
  W^s_\epsilon(p)\right)\leq 2n+\frac{b^+_i}{s^{+}_i}< 2n
\end{equation}
for all $i=1,\dots,m$. This implies
\begin{equation}
\overline{\dim}_B\left(\bigcup_{i=1}^m \bigcup_{p\in
J_i}W^s_\epsilon(p)\right)\leq 2n+\max\{b^+_i/s^+_i\}<2n.
\end{equation}
Hence
\begin{equation}
\dim_H\left(\bigcup_{p\in
J}W^s_\epsilon(p)\right)\leq 2n+\max\{b^+_i/s^+_i\}<2n.
\end{equation}
Since
\begin{equation}
W^s(p)=\bigcup_{k\in \bN}f^{-k}(W^s_\epsilon(f^k(p)))
\end{equation}
for all $p\in J$, inequality \eqref{eqghj} follows from Theorem
\ref{thghj} (i), Proposition \ref{l:hyper} (iii) and the fact that the Hausdorff dimension is countable
union stable.
\end{proof}
\noindent {\it Remark. } In the case $n=2$ the result $\dim_H
J^\pm<4$ was already shown in \cite{W1}.  However, the methods used in
\cite{W1} crucially depend on the fact that the unstable/stable
index of $J$ is identically one, and therefore,  do not apply 
to the case $n>2$.

\end{document}